\documentclass[11pt,reqno]{amsart}

\usepackage[english]{babel}
\usepackage{amsmath,amssymb,amsthm}
\usepackage{bbm}
\usepackage{calrsfs}
\usepackage{mathrsfs}

\usepackage[pdftex]{color}
\usepackage[dvipsnames]{xcolor}
\usepackage[bookmarks=true,hyperindex,pdftex,colorlinks,
citecolor=Black,linkcolor=Black, urlcolor=Black
]{hyperref}

\usepackage{ulem} 

\numberwithin{equation}{section}

\theoremstyle{plain}
\newtheorem{theorem}{Theorem}[section]
\newtheorem{theo}[theorem]{Theorem}

\newtheorem{prop}[theorem]{Proposition}
\newtheorem{lemma}[theorem]{Lemma}

\theoremstyle{definition}

\newtheorem{example}[theorem]{Example}

\newtheorem{defi}[theorem]{Definition}

\newtheorem{problem}[theorem]{Problem}

 \DeclareMathOperator{\ex}{ex}

\newcommand{\R}{\mathbb{R}}
\newcommand{\N}{\mathbb{N}}

\newcommand{\eins}{\mathbbm{1}}

\newcommand{\dopu}{{:}\allowbreak\ }
\newcommand{\eps}{\varepsilon}
\newcommand{\Id}{\mathrm{Id}}

\def\bea{\begin{align*}}
\def\beq{\begin{equation}}
\def\eeq{\end{equation}}
\newcommand{\begsta}{\begin{statements}}
\def\endsta{\end{statements}}

\newcounter{abc}

\newenvironment{statements}%
{\setcounter{abc}{0}
\begin{list}%
{{\rm (\alph{abc})}}
{\usecounter{abc}
\parsep=0pt plus 1pt
\topsep=1pt plus 2pt minus 1pt
\itemsep=1pt plus 2pt minus 1pt
\leftmargin=3\baselineskip
\labelsep=.6\baselineskip
\labelwidth=2.4\baselineskip
\rightmargin 0pt}%
}%
{\end{list}}

\makeatletter
\@namedef{subjclassname@2020}{%
  \textup{2020} Mathematics Subject Classification}
\makeatother

\begin{document}

\title[On Wolfgang Lusky's paper ``The Gurarij spaces are unique'']{On Wolfgang Lusky's paper\\``The Gurarij spaces are unique''}

\author{Dirk Werner}

\address{Department of Mathematics \\ Freie Universit\"at Berlin \\
Arnimallee~6 \\ D-14195~Berlin \\ Germany\newline
\href{http://orcid.org/0000-0003-0386-9652}{ORCID: \texttt{0000-0003-0386-9652}}
}
\email{werner@math.fu-berlin.de}



\begin{abstract}
This note surveys Wolfgang Lusky's proof of uniqueness of the Gurariy spaces and mentions further developments. 
\end{abstract}

\date{7.8.23; version v2}

\thanks{This piece has been commissioned by the editors of Archiv der Mathematik on the occasion of the 75th anniversary of the journal}

\maketitle
\thispagestyle{empty} 


\section{Introduction}\label{sec1}

In 1966, V.~I.~Gurariy \cite{Gu} defined the notion of a \textit{Banach space of (almost) universal disposition} by a certain extension property; see Definition~\ref{def2.1}. He proved the existence of (separable) such spaces and investigated some of their properties; henceforth, such spaces were called \textit{Gurariy spaces} (alternative spellings: Gurarii, Gurarij, Gurari\u \i, \dots); we shall reserve this name to separable spaces of this kind. While it is not a daunting task to prove that any two Gurariy spaces are almost isometric in the sense that their Banach-Mazur distance is~$1$, it remained open to decide whether they are actually isometric. This was asked for instance by J.~Lindenstrauss and his collaborators at various junctures (\cite[Problem~II.4.13]{LiTz0}, \cite{LazLin}).

The isometry problem was solved in 1976 by a fresh PhD from the (likewise rather freshly established) University of Paderborn, Wolfgang Lusky, in his first-ever published paper (the title says it all)
\begsta
\item[{[L]}]
The Gurarij spaces are unique. \textit{Arch.\ Math.}\ 27, 627--635 (1976).
\endsta
We shall refer to this paper, which is \cite{L} in the bibliography, simply by~[L].

The present note aims at surveying the background, Lusky's proof, and the ramifications of this result along with an outlook. 

Interestingly, some 30~years later Gurariy and Lusky cooperated intensively on a rather different topic, the M\"untz spaces, which has led to their monograph \cite{GuL}.

The notation in this note is standard; $B_X$ stands for the closed unit ball of $X$ and $\ex B_X$ for the set of its extreme points. We are considering only real Banach spaces. 


\section{Banach spaces of almost universal disposition}\label{sec2}

V.~I.~Gurariy (1935--2005)  was a member of the Kharkiv school of Banach spaces led by M.~I.~Kadets (sometimes spelled Kadec), one of the strongest in Europe which had its heyday from the late 1950ies till the collapse of the Soviet Union that produced a brain-drain in all fields of science. Gurariy himself emigrated to the United States in the early 1990ies.  After 2000, the Kharkiv school was basically reduced to V.~Kadets and his students. In 2022 the terror regime in Moscow set out to destroy the university of Kharkiv altogether \cite{Notices}, but remembering a slogan from many years back, !`No pasar\'an! 

Here is the key definition of his paper \cite{Gu}.

\begin{defi}\label{def2.1}
Let $X$ be a Banach space with the following property.
\begsta
\item[$\bullet$]
For finite-dimensional spaces $E$ and $F$, isometries $T\dopu E\to X$ and $S\dopu E\to F$, and for $\eps>0$, there exists an operator $\widehat T\dopu F\to X$  satisfying $\widehat TS=T$ and
\[
(1+\eps)^{-1} \|y\| \le \|\widehat Ty\| \le (1+\eps)\|y\| \qquad  (y\in F)
\]
(``an $\eps$-isometry''). 
\endsta
Then $X$ is called a Banach space of \textit{almost universal disposition}. A separable such space will also be called a \textit{Gurariy space}. 
\end{defi}

The epithet ``almost'' in this definition refers to the quantifier ``for all $\eps>0$''; if $\eps=0$ is permissible above, then the ``almost'' will be dropped. However, Gurariy proved in \cite[Th.~10]{Gu} that no separable space of universal disposition exists, but see Subsection~\ref{subsec6.3}  below.

If in the above definition, $S$ is the identical inclusion, i.e., $E\subset F$, then $\widehat T$ is an extension of $T$, which can likewise be considered as the identical inclusion. 

To see that the condition of Definition~\ref{def2.1} is quite restrictive, let us discuss two examples.

\begin{example}
(a) $c_0$ is not a space of almost universal disposition. Indeed, let $E=\R$, $T\dopu E\to c_0$, $T(r)=(r,0,0,\dots)$, $F=\ell_\infty^2 = \R^2$ with the max-norm, $S\dopu E\to F$, $S(r)=(r,r)$. Assume that $\widehat T$ has the properties of Definition~\ref{def2.1}, and let $\widehat T(-1,1)= (x_1, x_2, \dots)$. Note that $\widehat T(1,1)= (1,0,0,\dots)$ and therefore
\begin{align*}
\widehat T(0,1) &= \Bigl( \frac{1+x_1}2, \frac{x_2}2, \dots \Bigr) , \\
\widehat T(1,0) &= \Bigl( \frac{1-x_1}2, \frac{-x_2}2, \dots \Bigr) .
\end{align*}
This shows that $\widehat T$ cannot be an $\eps$-isometry for small $\eps$. 
(If $x$ is a real number close to~$1$ in modulus, then $\frac{1\pm x}2$ cannot both be close to~$1$.)

(b)  $C[0,1]$ is not a space of almost universal disposition. Indeed, let $E=\R$, $T\dopu E\to C[0,1]$, $T(r) = r\eins$ (the constant function), $F= \ell_2^2 = \R^2$ with the $\ell_2$-norm, $S\dopu E\to F$, $S(r)=(r,0)$. Assume that $\widehat T$ has the properties of Definition~\ref{def2.1}, and let $\widehat T(0,1)= f$. Note that $\widehat T(1,0)= \eins$ and therefore
$$
\widehat T(1,1)= \frac{\eins + f}2,
$$
which must have norm $\sqrt2 = \|(1,1)\|_2$ up to $\eps$. Since $(1+\eps)^{-1} \le \|f\| \le 1+\eps$, this is impossible for small~$\eps$. 
\end{example}

These examples indicate that positive results might not be very easy to come by. By a technical inductive argument, Gurariy shows in \cite[Th.~2]{Gu} the following existence theorem. 

\begin{theo}\label{th2.3}
There exists a separable Banach space of almost universal disposition. 
\end{theo}

As for uniqueness, he proves the following result. To formulate it succinctly, let us recall the \textit{Banach-Mazur distance} between (isomorphic) Banach spaces
$$
d(X,Y)= \inf\{ \|\Phi\| \|\Phi^{-1}\|\dopu \ \Phi\dopu X\to Y \text{ is an isomorphism}\}
$$
and call two Banach spaces \textit{almost isometric} if their Banach-Mazur distance equals~$1$.

Now for Theorem~5  of \cite{Gu}.

\begin{theo}\label{th2.4}
Any two separable spaces of almost universal disposition are almost isometric.
\end{theo}

A quick sketch of the proof can also be found in \cite[p.~168]{LiTz0}.


\section{The Lazar-Lindenstrauss approach}\label{sec3}

A key property of the Gurariy spaces (from now on we shall use this terminology) is that they are $L_1$-preduals. Recall that an \textit{$L_1$-predual} (a.k.a.\ a \textit{Lindenstrauss space}) is a Banach space whose dual is isometrically isomorphic to a space $L_1(\mu)$ of integrable functions on some measure space. This class of spaces is the subject of Lindenstrauss's epoch-making memoir \cite{Lin-mem}. 

\begin{prop}\label{prop3.1}
Every Gurariy space is an $L_1$-predual.
\end{prop}

In the literature, especially from the previous century, there are only vague indications as to why this is so. Since a recent article \cite{CFR}  admits that this proposition is ``not completely evident from the definition'' and since it is  instrumental for Lusky's proof, I'll sketch a proof. To begin with, we have to recall a characterisation of $L_1$-preduals from Lindenstrauss's memoir; see \cite[Th.~6.1]{LiTz0} in conjunction with \cite[Lemma~4.2]{LiTz0}, or \cite[\S21]{Lacey}.

\begin{theo}\label{th3.2}
A Banach space $X$ is an $L_1$-predual if and only if any four open balls $U(x_i, r_i)$ that intersect pairwise have a nonvoid intersection. It is enough to check this for balls of radius~$1$. 
\end{theo}

Let us verify that a Gurariy space $X$ has this property. So suppose $U(x_1,1), \dots , U(x_4,1)$ are four open balls in $X$ with radius~$1$ that intersect pairwise, i.e., $\|x_i-x_j\|<2$. Choose $\eps>0$ such that even $\|x_i-x_j\|<2-4\eps$. Let $E$ be the span of $x_1,\dots,x_4$. There are some $N\in \N$ and a linear operator $S_1\dopu E\to \ell_\infty^N$ such that 
$$
\frac1{1+\eps} \|S_1x\|_\infty \le \|x\| \le \|S_1x\|_\infty \qquad (x\in E).
$$
Let us consider the balls $U_{\ell_\infty^N}(S_1x_i, 1-\eps)$ in $\ell_\infty^N$. They intersect pairwise since
$$
\|S_1x_i-S_1x_j\|_\infty \le (1+\eps) \|x_i-x_j\| < (1+\eps)(2-4\eps) < 2-2\eps.
$$
Being pairwise intersecting balls in $\ell_\infty^N$, these balls have a point in common. This means that there exists some $z\in \ell_\infty^N$ such that 
$$
\|z-Sx_i\|_\infty <1-\eps \qquad (i=1,\dots,4).
$$
Unfortunately, $S_1$ is not an isometry and therefore is not eligible for being used in Definition~\ref{def2.1}. However, we can renorm $\ell_\infty^N$ to make it an isometry: note that $B_{\ell_\infty^N} \cap S_1(E) \subset S_1(B_E)$, and we can renorm $\ell_\infty^N$ by letting the new unit ball be the convex hull of $ S_1(B_E)$ and $B_{\ell_\infty^N}$. Call this renorming $F$, and let $S=S_1$ considered as an operator from $E$ to $F$; this is an isometry. We have 
$$
\frac1{1+\eps} \|y\|_\infty \le \|y\|_F \le \|y\|_\infty \qquad (y\in F)
$$
and thus 
$$
\|z-Sx_i\|_F \le \|z-S_1x_i\|_\infty < 1-\eps.
$$
Since $X$ is a Gurariy space, there is an $\eps$-isometry $\widehat T\dopu F\to X$ satisfying $\widehat TSx=x$ for $x\in E$. Let $x_0=\widehat Tz$; then $x_0\in \bigcap U(x_i,1)$:
$$
\|x_0-x_i\| = \|\widehat T z - \widehat T Sx_i\| \le (1+\eps) \|z-Sx_i\|_F <(1+\eps)(1-\eps) <1.
$$

In the more contemporary literature one can find explicit proofs of Proposition~\ref{prop3.1} based on another characterisation of $L_1$-preduals and a ``pushout argument'' \cite[Th.~2.17]{GarKub}, \cite[Prop.~6.2.8]{CC}.

Now let $X$ be a separable $L_1$-predual. By the results of Michael and Pe\l czy\'nski \cite{MicPel} and Lazar and Lindenstrauss \cite{LazLin} there is a chain of finite-dimensional subspaces $E_n$ of $X$ such that 
\begsta
\item
$E_1\subset E_2 \subset \dots$;
\item
$\dim E_n=n$, and $E_n$ is isometrically isomorphic to $\ell_\infty^n$,
\item
$\bigcup E_n$ is dense in $X$.
\endsta
The inclusion $E_n\subset E_{n+1}$ entails some degree of freedom, namely the choice of an isometry $\psi_n\dopu \ell_\infty^n \to \ell_\infty^{n+1}$. To study the structure of these $\psi_n$, we need the ad-hoc notion of an admissible basis: if $\delta_1,\dots,\delta_n$ denotes the canonical unit vector basis of $\ell_\infty^n$ and $\psi\dopu \ell_\infty^n \to \ell_\infty^n$ is an isometry, then $\psi(\delta_1),\dots,\psi(\delta_n)$ is called an \textit{admissible basis} for $\ell_\infty^n$. Note that $\psi$ takes a vector $(a_1,\dots,a_n)$ to $(\vartheta_1 a_{\pi(1)}, \dots, \vartheta_n a_{\pi(n)})$ for some permutation $\pi$ and some signs $\vartheta_j=\pm1$. Thus, an admissible basis  is just a permutation of the unit vector basis up to signs, and the isometric image of an admissible basis is again an admissible basis. 

Let us return to the isometric embedding $\psi_n\dopu \ell_\infty^n \to \ell_\infty^{n+1}$, and let $e_{1,n}, \dots , \allowbreak e_{n,n}$ be an admissible basis for~$\ell_\infty^n$. We can develop the vectors $f_j:= \psi_n(e_{j,n})$ into the unit vector basis of~$\ell_\infty^{n+1}$. Since $\psi_n$ is an isometry, there is at least one coordinate~$i$ where $|f_j(i)|=1$. Then, if $k\neq j$, $f_k(i)=0$: pick a sign $\lambda$ such that 
$$
|f_j(i) +\lambda f_k(i)| = |f_j(i)| + |f_k(i)| = 1+ |f_k(i)| 
$$
and so
$$
1= \|e_{j,n} + \lambda e_{k,n}\| = \|f_{j} + \lambda f_{k}\| \ge |f_{j}(i) + \lambda f_{k}(i) | =
1+ |f_k(i)|,
$$
hence the claim. Since $\|\psi_n\|=1$, we also have
$$
\Bigl| \sum_{j=1}^n f_j(i) \Bigr| = \Bigl| \psi_n \Bigl( \sum_{j=1}^n e_{j,n} \Bigr) (i) \Bigr| \le
\Bigl\| \sum_{j=1}^n e_{j,n} \Bigr\| =1.
$$
Therefore, there is an admissible basis $e_{1,n+1},  \dots, e_{n+1, n+1}$ for $\ell_\infty^{n+1}$ such that for some numbers $a_{jn}$
$$
\psi_n(e_{j,n}) = e_{j,n+1} +a_{jn} e_{n+1, n+1} \qquad (j=1,\dots, n)
$$
and 
$$
\sum_{j=1}^n |a_{jn}| \le 1.
$$

We can rephrase these representations in terms of the $E_n$ as follows. 

\begin{prop}\label{prop3.3}
There exist admissible bases in each $E_n$ and real numbers $a_{jn}$ such that 
\begin{align*}
&e_{j,n} = e_{j,n+1} +a_{jn} e_{n+1, n+1} \qquad (j=1,\dots, n; \ n=1,2,\dots) \\
&\sum_{j=1}^n |a_{jn}| \le 1 \qquad (n=1,2,\dots).
\end{align*}
\end{prop}

This proposition is due to Lazar and Lindenstrauss \cite{LazLin}.
The triangular matrix $(a_{jn})_{j\le n, n\in \N}$ is called a \textit{representing matrix} for the given $L_1$-predual~$X$. Conversely does the choice of admissible bases and of an array $(a_{jn})$ lead to an $L_1$-predual.

 Lazar and Lindenstrauss use this approach to present another proof of the existence of Gurariy spaces. Let $a_n=(a_{1n}, \dots , a_{nn}, 0,0,\dots)$ be the $n^{\mathrm{th}}$ column of a matrix as in Proposition~\ref{prop3.3}; then each $a_n$ is in the unit ball of~$\ell_1$. 

\begin{theo}\label{th3.4}
If $\{a_1, a_2, \dots\}$ is dense in the unit ball of $\ell_1$, then the corresponding matrix is associated to a Gurariy space. 
\end{theo}

It should be noted that the representing matrix $A$ of an $L_1$-predual $X$ is not uniquely determined, and much work has been done  to study the relation of $A$ and $X$ for certain classes of $L_1$-preduals; see e.g.\ Lusky's paper \cite{L77}.


\section{Lusky's uniqueness proof}\label{sec4}

Here is Lusky's uniqueness theorem.

\begin{theo}\label{th4.1}
Any two Gurariy spaces are isometrically isomorphic.
\end{theo}

Let us first remark that almost isometric spaces (cf.\ Theorem~\ref{th2.4}) need not be isometric. The following is a classical counterexample due to  Pe\l czy\'nski  from \cite{Pel}: Let $X$ and $Y$ be $c_0$ equipped with the equivalent norms ($x=(x_n)$)
\begin{align*}
\|x\|_X &= \|x\|_\infty + \Bigl( \sum_{n=1}^\infty \frac{|x_{n}|^2}{2^n} \Bigr)^{1/2} , \\
\|x\|_Y &= \|x\|_\infty + \Bigl( \sum_{n=1}^\infty \frac{|x_{n+1}|^2}{2^n} \Bigr)^{1/2}.
\end{align*} 
The operators $\Phi_n\dopu X\to Y$, $x\mapsto (x_n, x_1, \dots, x_{n-1}, x_{n+1}, \dots)$ are isomorphisms satisfying $\lim_n \|\Phi_n\| \|\Phi_n^{-1}\| =1$ so that $X$ and $Y$ are almost isometric; but $X$ is strictly convex while $Y$ isn't, therefore $X$ and $Y$ are not isometric. 

Benyamini \cite{Ben} has shown that such counterexamples also exist among $L_1$-preduals.

The proof of Theorem~\ref{th4.1} consists of a delicate inductive construction of $\ell_\infty^n$-subspaces and admissible bases.  The key problem to be solved here is this. 

\begin{problem}\label{prob4.2}
Let $X$ be a Gurariy space and $E\subset F$ be finite-dimensional spaces with $E\cong \ell_\infty^n$ and $F\cong \ell_\infty^{n+1}$. Let $T\dopu E\to X$ be an isometry. When does there exist an isometric extension $\widehat T\dopu F\to X$?
\end{problem}

Lusky notes that this is not always the case [L, p.~630], and he gives the following useful criterion in terms of admissible bases. W.l.o.g.\ suppose that $T$ is the identity. Let $e_1,\dots, e_n$ and $f_1,\dots,f_{n+1}$ be admissible bases for $E$ resp.\ $F$ such that 
$$
e_i = f_i + r_i f_{n+1}, \qquad i=1,\dots,n.
$$

\begin{lemma}\label{lem4.3}
Problem~\ref{prob4.2} has a positive solution if $\sum_{i=1}^n |r_i| < 1$. 
\end{lemma}

This criterion is a little hidden in the proof of the Corollary [L, p.~630], where the extreme point condition $\ex B_E \cap \ex B_F= \emptyset$ is spelled out to be sufficient; but the heart of the matter is Lemma~\ref{lem4.3}.

Now let's take a quick glimpse at the proof of Theorem~\ref{th4.1}. Suppose that $X$ and $Y$ are Gurariy spaces coming with $\ell_\infty^n$-approximations $\bigcup_n E_n$ and $\bigcup F_n$, respectively. Comparing Proposition~\ref{prop3.3} with Lemma~\ref{lem4.3} one realises that one has to perturb the given admissible bases so that Lemma~\ref{lem4.3} becomes applicable. The details of this process are quite technical [L, pp.~631--633] and lead to sequences of admissible bases. Ultimately one can pass to the limit and obtain admissible bases $\{e_{i,n}\dopu i\le n$, $n\ge1\}$ resp.\ $\{f_{i,n}\dopu i\le n$, $n\ge1\}$ spanning dense subspaces of $X$ resp.\ $Y$, and the operator $e_{i,n} \mapsto f_{i,n}$ acts as a well-defined isometry. 

In an addendum to [L], dated January 10, 1976, Lusky applies his methods to Mazur's rotation problem that asks  whether a separable transitive space is isometric to a Hilbert space; a Banach space $X$ is called \textit{transitive} if whenever $\|x\|=\|y\|=1$, there is an isometric automorphism $T\dopu X\to X$ mapping $x$ to $y$, i.e., $Tx=y$. This problem is open to this day, and recent papers on the subject include \cite{Cab} and \cite{CFR}.

What Lusky proves in his addendum  is that the Gurariy space (now that we know it's unique we may use the definite article) is transitive for smooth points. Recall that $x_0$ is a smooth point of the unit ball $B_X$ if $\|x_0\|=1$ and there is exactly one $x_0^*\in X^*$ such that $\|x_0^*\| = x_0^*(x_0)=1$; equivalently, the norm function $x\to \|x\|$ is G\^ateaux differentiable at~$x_0$. It is a theorem of Mazur that smooth points are dense in the unit sphere of a separable Banach space. 

\begin{theo}\label{th4.4}
Let $x$ and $y$ be smooth points of the unit ball of the Gurariy space~$G$. Then there is  an isometric automorphism $T\dopu G\to G$ mapping $x$ to $y$.
\end{theo}

Another result of [L] is a refined version of a theorem originally due to Wojtaszczyk \cite{Woj} (see also \cite{L77}).

\begin{theo}\label{th4.5}
Let $X$ be a separable $L_1$-predual and $G$ be the Gurariy space. Then there exist an isometry $T\dopu X\to G$ and a norm-$1$ projection $P\dopu G\to G$ onto $T(X)$; further $(\Id-P)(G)$ is isometrically isomorphic to~$G$.
\end{theo}

This indicates that the Gurariy space is ``maximal'' among the separable $L_1$-predual spaces; in particular it contains $C[0,1]$ and is universal, a fact proved by other means by Gevorkyan in \cite{Gev}. 

We close this section by mentioning another proof of Theorem~\ref{th4.1}, due to W.~Kubi\'s and S.~Solecki \cite{KubSol}. Their proof avoids the Lazar-Lindenstrauss machinery and just depends on the defining properties of a Gurariy space. They also prove the universality of the Gurariy space from first principles, without relying on the universality of $C[0,1]$. Still another proof is in Kubi\'s's paper \cite{Kub} in \textit{Archiv der Mathematik}, which builds on a Banach-Mazur type game.


\section{The Poulsen simplex}\label{sec5}

This note wouldn't be complete without mentioning the cousin of the Gurariy space in the world of compact convex sets, the \textit{Poulsen simplex}. The traditional definition of a (compact) simplex is a compact convex subset $S$ of a Hausdorff locally convex space $E$ such that the cone generated by $S\times \{1\}$ in $E\oplus \R$ is a lattice cone. Thus, a triangle in the plane is a simplex while a rectangle isn't. For our purposes it is important to note that the space $A(S)$ of affine continuous functions on a compact convex set is an $L_1$-predual if and only if $S$ is a simplex. 

Poulsen \cite{Pou} had proved the existence of a metrisable simplex, which now bears his name, whose set of extreme points is dense. It is a result due to Lindenstrauss, Olsen, and Sternfeld \cite{LOS} that such a simplex is uniquely determined up to affine homeomorphism. They write:  

\begin{quote}\small
We discovered the uniqueness of the Poulsen simplex after reading Lusky's paper [L] on the uniqueness of the Gurari space. Our proof of the uniqueness uses the same idea which Lusky used in~[L].
\end{quote}
The role of admissible bases is now played by peaked partitions of unity. 

The authors mention a lot of similarities between the Poulsen simplex and the Gurariy space. For example, the counterpart of the defining property of the Poulsen simplex $S_P$ is Lusky's theorem from [L] and \cite{L77} that a separable $L_1$-predual is a Gurariy space $G$ if and only if $\ex B_{G^*}$ is weak$^*$ dense in the unit ball $B_{G^*}$. However, $A(S_P)$ is not the Gurariy space since for example the transitivity property of Theorem~\ref{th4.4} fails. But, as shown by Lusky \cite{L-LOS}, one can salvage this by requiring a slightly more stringent assumption on $x$ and $y$, which are now supposed to be positive: in addition, $\eins-x$ and $\eins-y$ should be smooth points.


\section{Outlook}\label{sec6}

\subsection{Fra\"\i ss\'e theory}\label{subsec6.1}
The Gurariy space is a very homogeneous object, for example \cite[Th.~3]{Gu}: If $E$ and $F$ are finite-dimensional subspaces of the same dimension of a Gurariy space $G$, then for every $\eps>0$, every isometric isomorphism from $E$ to $F$ extends to an $\eps$-isometric automorphism of~$G$. In recent  years, such homogeneous structures were investigated by methods of model theory known as Fra\"\i ss\'e theory (\cite{Fra}, \cite{BenYaa}, \cite{Kub14}). Fra\"\i ss\'e theory associates a unique limit to certain substructures. This approach is at least implicit in the Kubi\'s-Solecki uniqueness proof, and a detailed exposition involving the Gurariy space, the Poulsen simplex and a whole lot more can be found in M.~Lupini's paper~\cite{Lup18}.

\subsection{Noncommutative Gurariy spaces}\label{subsec6.2}
T.~Oikhberg, in his \textit{Archiv der Mathematik} paper \cite{Oik}, proved the existence and uniqueness of a ``noncommutative'' Gurariy space, i.e., a Gurariy-like object in the setting of operator spaces \`a la Effros-Ruan.  Again, this can also be viewed from the perspective of Fra\"\i ss\'e theory~\cite{Lup16}.

\subsection{Nonseparable spaces}\label{subsec6.3}
We have already mentioned in Section~\ref{sec2} Gurariy's result that no space of universal disposition can be separable. Since the definition of (almost) universal disposition makes perfect sense beyond the separable case, it was studied in several papers, e.g., \cite{Avi-etc}, \cite{CasSim}, \cite{GarKub}. It turns out that there are spaces of almost universal disposition of density character $\aleph_1$, but the uniqueness breaks down (Th.~3.6 and Th.~3.7 in \cite{GarKub}). Likewise, there are spaces of universal disposition of density $\aleph_1$, and again, uniqueness fails (\cite{Avi-etc}, \cite{CasSim}). Indeed, it should be noted that in these papers also the variant of being of (almost) universal disposition with respect to separable spaces, already considered by Gurariy, is studied: in Definition~\ref{def2.1} one now allows $E$ and $F$ to be separable rather than finite-dimensional.

\subsection{Banach lattices}\label{subsec6.4}
Recently, M.~A.~Tursi \cite{Tur} proved the existence of a uniquely determined Gurariy-like Banach lattice. She exploits ideas of Fra\"\i ss\'e theory.



\begin{thebibliography}{99}


\bibitem{Avi-etc}
\textsc{A.~Avil\'es}, \textsc{F.~Cabello S\'anchez},  \textsc{J.~M.~F. Castillo},  \textsc{M.~Gonz\'alez}, and  \textsc{Y.~Moreno}, 
Banach spaces of universal disposition.
\emph{J. Funct. Anal.} 261, No.~9 (2011), 2347--2361. 
Zbl 1236.46014


\bibitem{BenYaa}
\textsc{I. Ben Yaacov},
Fra\"\i ss\'e limits of metric structures. 
\emph{J. Symb. Log.} 80, No.~1 (2015), 100--115. 
Zbl 1372.03070



\bibitem{Ben}
\textsc{Y.~Benyamini},
Near isometries in the class of $L^1$-preduals.
\emph{Isr. J. Math.} 20 (1975), 275--281. 
Zbl 0314.46016


\bibitem{Cab}
\textsc{F. Cabello S\'anchez},
Wheeling around Mazur rotations problem. 
\emph{Trans. Amer. Math. Soc.} 376, No.~3 (2023), 2213--2235. 
Zbl 07655332


\bibitem{CC}
\textsc{F. Cabello S\'anchez} and \textsc{J.~M.~F. Castillo},
Homological Methods in Banach Space Theory. 
Cambridge Studies in Advanced Mathematics 203. Cambridge: Cambridge University Press 2023. 
Zbl 07625516


\bibitem{CFR}
\textsc{F. Cabello S\'anchez}, \textsc{V.~Ferenczi}, and \textsc{B.~Randrianantoanina},
On Mazur rotations problem and its multidimensional versions. 
\emph{S\~ao Paulo J. Math. Sci.} 16, No.~1 (2022), 406--458. 
Zbl 1502.46004


\bibitem{CasSim}
\textsc{J.~M.~F. Castillo} and \textsc{M.~Sim\~oes},
On Banach spaces of universal disposition. 
\emph{New York J. Math.} 22 (2016), 605--613. 
Zbl 1355.46009


\bibitem{Fra}
\textsc{R. Fra\"\i ss\'e},
Sur l’extension aux r\'elations de quelques propri\'et\'es des ordres. 
\emph{Ann. Sci. \'Ec. Norm. Sup\'er.}, III.~S\'er. 71 (1954), 363--388.
Zbl 0057.04206


\bibitem{GarKub}
\textsc{J.~Garbuli\'nska} and \textsc{W.~Kubi\'s},
Remarks on Gurari\u\i\ spaces. 
\emph{Extr. Math.} 26, No.~2 (2011), 235--269. 
Zbl 1267.46020


\bibitem{Gev}
\textsc{Yu. L. Gevorkyan},
The universality of spaces of almost universal placement.
\emph{Funct. Anal. Appl.} 8 (1974), 157; translation from \emph{Funkts. Anal. Prilozh.} 8, No.~2 (1974), 72. 
 Zbl 0296.46019


\bibitem{Gu}
\textsc{V. I. Gurari\u\i},
Space of universal disposition, isotropic spaces and the Mazur problem on rotations of Banach spaces. 
\emph{Sib. Math. J.} 7 (1966), 799--807 (1967); translation from \emph{Sib. Mat. Zh.} 7 (1966), 1002--1013.
Zbl 0166.39303 


\bibitem{GuL}
\textsc{V.~I. Gurariy} and \textsc{W.~Lusky},
Geometry of M\"untz spaces and related questions. 
Lecture Notes in Mathematics 1870. Springer Berlin 2005. 
Zbl 1094.46003


\bibitem{Kub14}
\textsc{W.~Kubi\'s},
Fra\"\i ss\'e sequences: category-theoretic approach to universal homogeneous structures.
\emph{Ann. Pure Appl. Logic} 165, No.~11 (2014), 1755--1811. 
 Zbl 1329.18002

\bibitem{Kub}
\textsc{W.~Kubi\'s},
Game-theoretic characterization of the Gurarii space. 
\emph{Arch. Math.} 110, No.~1 (2018), 53--59. 
Zbl 1396.46007


\bibitem{KubSol}
\textsc{W.~Kubi\'s} and \textsc{S.~Solecki},
A proof of uniqueness of the Gurari\u\i\ space.
\emph{Isr. J. Math.} 195, Part~A (2013), 449--456. 
Zbl 1290.46010


\bibitem{Lacey}
\textsc{H. E. Lacey},
The Isometric Theory of Classical Banach Spaces. Grundlehren Vol.~208. Springer Berlin-Heidelberg-New York 1974.
Zbl 0285.46024


\bibitem{LazLin}
\textsc{A. Lazar} and \textsc{J.~Lindenstrauss},
Banach spaces whose duals are $L_1$ spaces and their representing matrices. 
\emph{Acta Math.} 126 (1971), 165--193. 
Zbl 0209.43201


\bibitem{Lin-mem}
\textsc{J.~Lindenstrauss},
Extension of compact operators.
\emph{Mem. Amer. Math. Soc.} 48  (1964). 
Zbl 0141.12001


\bibitem{LOS}
\textsc{J. Lindenstrauss}, \textsc{G. Olsen}, and \textsc{Y. Sternfeld},
The Poulsen simplex. 
\emph{Ann. Inst. Fourier} 28, No.~1 (1978), 91--114. 
Zbl 0363.46006


\bibitem{LiTz0}
\textsc{J. Lindenstrauss} and \textsc{L.~Tzafriri},
Classical Banach Spaces.
Lecture Notes in Mathematics 338. Springer Berlin-Heidelberg-New York 1973.
Zbl 0259.46011


\bibitem{Lup16}
\textsc{M. Lupini},
Uniqueness, universality, and homogeneity of the noncommutative Gurarij space. 
\emph{Adv. Math.} 298 (2016), 286--324. 
Zbl 1348.46058


\bibitem{Lup18}
\textsc{M. Lupini},
Fra\"\i ss\'e limits in functional analysis. 
\emph{Adv. Math.} 338 (2018), 93--174. 
Zbl 1405.46041


\bibitem{L}
\textsc{W.~Lusky},
The Gurarij spaces are unique. 
\emph{Arch. Math.} 27 (1976), 627--635. 
Zbl 0338.46023


\bibitem{L77}
\textsc{W.~Lusky},
On separable Lindenstrauss spaces.
\emph{J. Funct. Anal.} 26 (1977), 103--120. 
Zbl 0358.46016


\bibitem{L-LOS}
\textsc{W.~Lusky},
A note on the paper `The Poulsen simplex’ of Lindenstrauss, Olsen and Sternfeld. 
\emph{Ann. Inst. Fourier} 28, No.~2 (1978), 233--243. 
Zbl 0335.46034


\bibitem{MicPel}
\textsc{E. Michael} and \textsc{A.~Pe\l czy\'nski},
Separable Banach spaces which admit $l_\infty^n$ approximations. 
\emph{Isr. J. Math.} 4 (1966), 189--198. 
Zbl 0151.17602


\bibitem{Oik}
\textsc{T. Oikhberg},
The non-commutative Gurarii space. 
\emph{Arch. Math.} 86, No.~4 (2006), 356--364. 
Zbl 1119.46045

\bibitem{Pel}
\textsc{A.~Pe\l czy\'nski} with the collaboration of \textsc{C.~Bessaga},
Some aspects of the present theory of Banach spaces. 
In: S.~Banach, \OE uvres, Vol.~2, PWN Warszawa (1979), pp.~223--302.
Zbl 0407.01009


\bibitem{Pou}
\textsc{E. T. Poulsen},
A simplex with dense extreme points. 
\emph{Ann. Inst. Fourier} 11 (1961), 83--87. 
Zbl 0104.08402


\bibitem{Tur}
\textsc{M. A. Tursi},
A separable universal homogeneous Banach lattice. 
\emph{Int. Math. Res. Not.} 2023, No.~7 (2023), 5438--5472. 
Zbl 07673008


\bibitem{Notices}
\textsc{M. Vlasenko} and \textsc{E. Zelmanov},
Voices from the bombed universities of Ukraine.
\emph{Notices Amer. Math. Soc.} June/July 2023, 987--991.


\bibitem{Woj}
\textsc{P.~Wojtaszczyk},
Some remarks on the Gurarij space. 
\emph{Studia Math.} 41 (1972), 207--210. 
Zbl 0233.46024


\end{thebibliography}
\end{document}